\documentclass[11pt]{article}
\usepackage{latexsym}

\usepackage{amsmath,amsfonts,amssymb}
\newcommand\RR{{\Bbb R}}
\newcommand\CC{{\Bbb C}}

\newcommand\ZZ{{\Bbb Z}}

\def\d={\,:=\,}

\newcommand{\semdir}
{\rtimes}

\font\frakten=eufm10
\newfam\frakfam
\textfont\frakfam=\frakten

\newtheorem{thm}{Theorem}[section]
\newtheorem{lemma}[thm]{Lemma}
\newtheorem{cor}[thm]{Corollary}
\newtheorem{prop}[thm]{Proposition}
\newtheorem{Defn}[thm]{Definition}
\newtheorem{Ex}[thm]{Example}
\newtheorem{Rem}[thm]{Remark}
\newtheorem{Exs}[thm]{Examples}
\newtheorem{Rems}[thm]{Remarks}
\newtheorem{Defrem}[thm]{Definition and Remark}
\newtheorem{Remnt}[thm]{}
\newenvironment{defn}
 {\begin{Defn} \begin{rm}} {\end{rm} \hfill $\Box$ \end{Defn}}

\newenvironment{prf} {{\bf Proof.}}{\hfill $\Box$}

\begin{document}

\title{Plancherel transform criteria for Weyl-Heisenberg frames with integer
oversampling}
\author{Hartmut F\"uhr \thanks{email: fuehr@gsf.de},
 \\ GSF Research Center for Environment and Health \\
 Institute for Biomathematics and Biometry \\
 D-85764 Neuherberg}
\date{\today}
\maketitle
\begin{abstract}
 We investigate the relevance of admissibility criteria based
 on Plancherel measure for the characterization of tight Weyl-Heisenberg frames
 with integer oversampling.
 For this purpose we observe that functions giving rise to such 
 Weyl-Heisenberg frames are admissible with regard to the action of a suitably
 defined type-I discrete group $G$. This allows to relate the construction of Weyl-Heisenberg frames
 to the Plancherel measure of $G$, which provides an alternative
 proof and a new interpretation of the well-known Zak transform based
 criterion for tight Weyl-Heisenberg frames with integer oversampling.
\end{abstract}

\section{Admissibility conditions and Weyl-Heisenberg frames}
\label{sec:Intro}

 This paper interprets characterizations of tight Weyl-Heisenberg frames as admissibility
 conditions connected to a certain discrete group. The starting point
 was an observed similarity between Zak transform based criteria for such 
 frames and representation-theoretic admissibility conditions established by the author.
 We will show that the former can be seen as special instances of the latter.

 In order to review the notion of admissibility, let $G$ be
 a locally compact group with left Haar measure $\mu_G$. We let
 ${\rm L}^2(G)$ denote the associated ${\rm L}^2$-space, on which $G$
 acts unitarily by left translations; this defines the {\bf left regular
 representation} $\lambda_G$ of $G$.
 Next let $(\pi,{\cal H}_\pi)$ be a strongly continuous unitary
 representation of $G$. Given vectors $\varphi, \eta \in {\cal H}_\pi$,
 the bounded continuous function 
 $V_{\eta} \varphi$ on $G$ is defined by $V_{\eta} \varphi(x) = \langle
 \varphi, \pi(x) \eta \rangle$. In the case that the linear
 mapping $V_\eta : \varphi \mapsto V_\eta \varphi$ defines an isometry
 ${\cal H}_\pi \hookrightarrow {\rm L}^2(G)$, we call $\eta$
 {\bf admissible}. The question which representations have admissible
 vectors has been answered in general for groups which have a well-behaved
 regular representation, i.e., for which $\lambda_G$ is type-I. The main
 device for proving these criteria is the Plancherel measure of the group,
 which also allows the characterization of admissible vectors. 
 The Plancherel measure $\nu_G$ underlies the decomposition
 of $\lambda_G$ into irreducibles; for a definition and background
 see \cite{Di,Fo}.

 The following theorem gives admissibility conditions for a restricted case,
 tailormade for the example which we study below. It is proved -- in greater
 generality -- in \cite{Fu}; the specialization to multiplicity-free representations
 is discussed in somewhat more detail in \cite{FuMa}.

 \begin{thm}
 \label{thm:adm_cond_Pl}
 Let $G$ by a unimodular locally compact group, and assume that
 $\lambda_G$ is type-I. Let $(\pi, {\cal H}_\pi)$ be a multiplicity-free
 representation. Then $\pi$ has admissible vectors iff there exists a 
 Borel subset $\Sigma \subset \widehat{G}$ with $\nu_G(\Sigma) < \infty$,
 such that
 \begin{equation} \label{eq:dir_int_dec}
 \pi \simeq \int_{\Sigma}^\oplus \sigma d\nu_G(\sigma) ~~.\end{equation} 
 If we identify $\pi$ with the right-hand side of (\ref{eq:dir_int_dec}),
 we obtain the following admissibility condition for vector fields
 $\eta = (\eta_\sigma)_{\sigma \in \Sigma} \in 
 {\cal H}_\pi$ :
 \begin{equation} \label{eqn:adm_cond_Pl}
 \eta \mbox{ is admissible } \Leftrightarrow \| \eta_\sigma \| = 1~~, \nu_G
 \mbox{-almost everywhere.}
 \end{equation}
 \end{thm}
 
 The theorem suggests the following three step program for establishing admissibility
 criteria for an arbitrary representation $\pi$:
\begin{itemize}
 \item Explicitly construct a unitary equivalence $T : {\cal H}_\pi
 \to \int_{\widehat{G}}^\oplus {\cal H}_\sigma d\widetilde{\nu}(\sigma)$, where
 $\widetilde{\nu}(\sigma)$ is a suitable measure on $\widehat{G}$.
 \item Compute $\nu_G$ and check whether $\widehat{\nu}$ is $\nu_G$-absolutely continuous.
 If the answer is no, or if $\widetilde{\nu}$ is not supported on a set $\Sigma$ of finite
 Plancherel measure, there is no admissible vector.
 \item If $\widetilde{\nu}$ is $\nu_G$-absolutely continuous and supported
 on a set of finite Plancherel measure, compute the Radon-Nikodym-derivative.
 This allows to compute
 the intertwining operator $\widetilde{T} : {\cal H}_\pi
 \to \int_{\Sigma}^\oplus {\cal H}_\sigma d\nu_G(\sigma)$ explicitly. Now
 $\eta \in {\cal H}_\pi$ is admissible iff $\widetilde{T} (\eta)
 \in \int_{\Sigma}^\oplus {\cal H}_\sigma d\nu_G(\sigma)$
 is a field of unit vectors. 
 \end{itemize}
 
 This approach may be aptly described as the exertion of representation-theoretic brute force:
 Both $\pi$ and $\lambda_G$ are decomposed into irreducibles, and then the containment of
 $\pi$ in $\lambda_G$ is checked by comparing the measures underlying the decompositions.
 The three steps require the computation of the unitary dual $\widehat{G}$, or at
 least of the support of the Plancherel measure. In addition, we need to explicitly
 compute a direct integral decomposition of $\pi$. 
 It is obvious that this may be quite inefficient when dealing
 with a concrete representation $\pi$, whose direct integral decomposition might be supported only on
 a small portion of $\widehat{G}$. As a matter of fact, this is what happens in the case we
 consider below. Moreover, in a concrete situation explicit knowledge of the 
 intertwining operators and measures may be hard to achieve.
 On the other hand, the approach provides a unified and complete 
 description of the representations, and puts the problem of devising
 admissibility conditions in a rather general representation-theoretic context.

 All technical problems notwithstanding, there are representations which have
 been analysed according to the general scheme.
 This includes a setting 
 which has been the object of a number of papers in wavelet theory, namely
 semidirect products of the type $\RR^k \semdir H$, where $H$ is a closed
 subgroup of ${\rm GL}(k,\RR)$. Such a group has a natural representation
 on ${\rm L}^2(\RR^k)$, the quasiregular representation arising from 
 the natural action on $\RR^k$. The problem of establishing admissibility
 conditions for these representations has been considered in varying degrees of
 generality in \cite{BeTa,Fu96,FuMa,LWWW}, and the connection to the Plancherel
 formula was worked out explicitly in \cite{FuMa}. Here the  
 Fourier transform of $\RR^k$ takes over the role of the
 intertwining operator $T$, and the computation of the
 measures is obtained by measure decomposition along the $H$-orbits on
 the character group $\widehat{\RR^k}$. In this paper, we 
 perform a similar analysis for the case of tight Weyl-Heisenberg frames with
 integer oversampling. It will turn out that in this case the Zak
 transform, which is the chief technical device in this setting,
 acts as the operator $T$. Moreover, the third step, 
 computing the Radon-Nikodym derivatives, turns out to be trivial.

 As in the case of the semidirect products, a disclaimer regarding the
 use of our scheme is in order. As was already mentioned above, applying the general results requires
 to check a fair amount of technical details, which do not necessarily contribute
 to a better understanding of the initial problem. By contrast, the usual proof of the Zak transform
 criterion is obtained by more or less elementary Fourier-analytic
 arguments. What we wish to stress is that the tight Weyl-Heisenberg frames considered in this paper
 constitute an illustration rather than an application of the abstract results. 

 In order to recall
 the definition of a Weyl-Heisenberg system, we define the translation 
 operators $T_x$ and modulation operators $M_\omega$ on ${\rm L}^2(\RR)$ by 
 \[ (T_x f) (y) = f(y-x) ~~,~~ (M_\omega f) (y) = e^{2 \pi i \omega y} f(y)
 ~~.\]
 Clearly these operators are unitary. Now a {\bf tight Weyl-Heisenberg frame}
 is a system of vectors
 \[ \psi_{i} = M_{\omega_i} T_{x_i} \psi ~~(i \in I),\]
 arising from a fixed vector $\psi \in {\rm L}^2(\RR)$ and a 
 countable family of phase-space points $(x_i,\omega_i) \in \RR^2$, ($i \in I$) such that for all
 $g \in {\rm L}^2(\RR)$, we have
 \begin{equation} \label{eqn:defn_ntf}
 \| g \|_2^2 = \sum_{i \in I} |\langle g, \psi_i \rangle|^2 ~~.\end{equation}
 There exist several alternative definitions, with varying indexing and
 ordering of operators. However, up to phase factors which clearly
 do not affect any of the frame properties, the resulting
 systems are identical.
 A {\bf normalized tight Weyl-Heisenberg frame with integer oversampling $\pmb{L}$}
 is a tight Weyl-Heisenberg
 frame arising from a family $\{ (x_{n,m},\omega_{n,m}) = (n,m/L); n,m \in \ZZ \}$.
 Now, given $L$, the problem is to decide for a given $\psi$ whether
 it induces a normalized tight Weyl-Heisenberg frame or not. As we will see in 
 the next section, the Zak transform allows a precise answer to this
 question. Our next aim is to show that the condition is in fact
 an admissibility condition for $\psi$. Note that for $L=1$,
 this is obvious: The set $\{ T_n M_m : n,m \in \ZZ \}$ is an
 abelian subgroup of the unitary group of ${\rm L}^2(\RR)$, and
 condition (\ref{eqn:defn_ntf}) precisely means admissibility in this case. For
 $L>1$ however,  $\{ T_n M_{m/L} : n,m \in \ZZ \}$ is not a subgroup,
 and we have to deal with the nonabelian group $G $ generated by this
 set. 

 From now on, we once and for all fix an integer oversampling rate $L \ge 1$,
 and define the underlying group $G$ as 
 \[ G  = \ZZ \times \ZZ \times (\ZZ / L \ZZ) ~~,\]
 with the group law
 \begin{equation} \label{eqn:group_law}
 (n,k,\overline{\ell}) (n',k',\overline{\ell'}) = (n+n',k+k',\overline{\ell}
 +\overline{\ell'}+\overline{k'n}) ~~\end{equation}
 and inverse given by $(n,k,\overline{\ell})^{-1}= (-n,-k,\overline{-\ell+kn})$.
 Here we used the notation $\overline{n} = n + L \ZZ$. The representation
 $\pi $ of $G $ acts on ${\rm L}^2(\RR)$ by 
 \[ \pi (n,k,\overline{\ell}) = e^{2 \pi i (\ell - nk)/L} M_{n/L} T_k = 
 e^{2 \pi i \ell/L} T_k M_{n/L}~~.\]
 It is straightforward to check how normalized tight frames with
 oversampling $L$ relate to admissibility for $\pi $:
 \begin{lemma} \label{lem:ntf_eq_admissible} Let
 $\psi \in {\rm L}^2(\RR)$. Then
 $(M_{n/L} T_k \psi)_{n,k \in \ZZ}$ is a normalized tight frame iff
  $\frac{1}{\sqrt{L}} \psi \in {\rm L}^2(\RR)$ is admissible for $\pi $.
 \end{lemma}
 \begin{prf}
 The relation
 \[ M_{n/L} T_k \psi = e^{-2 \pi i (\ell - nk)/L} \pi(k,n,\overline{\ell}) 
 \psi \]
 implies for all $g \in {\rm L}^2(\RR)$ that
 \[
  \sum_{n,k,\overline{\ell}} \left| \langle g, \pi(k,n,\overline{\ell}) f \rangle \right|^2
 = L \sum_{n,k} \left| \langle g, M_{n/L} T_k f \rangle \right|^2 ~~,
 \]
 which shows the claim.
 \end{prf}

 The following lemma establishes that $G $
 is a finite extension of an abelian normal subgroup $N$. It is important
 in two ways: It ensures that $G $ is type-I, and secondly, the
 computation of the dual $\widehat{G }$ and the Plancherel measure
 on it will be obtained by a Mackey analysis of this extension.

 \begin{lemma}
  Let 
 \[ N  = \{ ( nL, k, \overline{\ell}) : k,n,\ell \in \ZZ \} ~~.\]
 Then $N $ is an abelian normal subgroup of $G $ with
 $G  / N  \cong \ZZ / L \ZZ$. In particular, $G $ is type-I. 
 \end{lemma}
 \begin{prf}
 The statements concerning $N $ are obvious from (\ref{eqn:group_law});
 for the description of $G /N $ use the representatives
 $(0,0,0),(1,0,0),\ldots,(L-1,0,0)$ of the $N $-cosets.
 The type-I property of $G $ is immediate from this by 
 Thoma's theorem \cite{Th}.
 \end{prf} 

 One aspect which makes this example particularly attractive for the author
 is that it shows that the three step scheme is sometimes even applicable
 to discrete groups. One of the main shortcomings of the scheme is
 the restriction that the regular representation has to be type-I. While for
 many connected Lie groups (such as semisimple, nilpotent,...) the requirement
 is fulfilled, it is rather restrictive for discrete groups
 \cite{Ka}. 

\section{Zak transform criteria for tight Weyl-Heisenberg frames}

 In this section we introduce the Zak transform and formulate the
 criterion for normalized tight Weyl-Heisenberg frames. Our main reference for
 the following will be \cite{Gr}. 
 \begin{defn}
 For $f \in C_c(\RR)$, define the Zak transform of $f$ as the
 function ${\cal Z} f : \RR^2 \to \CC$ given by
 \[ {\cal Z} f (x,\omega) = \sum_{m \in \ZZ} f(x-m) e^{2 \pi i m \omega} ~~.
 \]
 \end{defn}
 
 The definition of the Zak transform immediately implies a quasi-periodicity
 condition for $F = {\cal Z} f$:
 \begin{equation} \label{eqn:per_cond}
 \forall m,n \in \ZZ~~:~~  F(x+m,\omega+n) = e^{2 \pi i m \omega} F(x,\omega) ~~.
 \end{equation}
 In particular, the Zak transform of a function $f$ is uniquely
 determined by its restriction to the unit square $[0,1]^2$. 
 We next extend the Zak transform to a unitary operator
 ${\cal Z}: {\rm L}^2(\RR) \to {\cal H}$, where ${\cal H}$ is a suitably
 defined Hilbert space. For the proof of the following see
 \cite[Theorem 8.2.3]{Gr}
 \begin{prop}
 Let the Hilbert space ${\cal H}$ be defined by
 \[ {\cal H} = \{ F \in {\rm L}^2_{\rm loc} (\RR^2 ) : F \mbox{ satisfies
 (\ref{eqn:per_cond}) almost everywhere on } \RR^2 \} ~~,
 \]
 with norm
 \[ \| F \|_{\cal H} = \| F \|_{{\rm L}^2([0,1]^2)} ~~.\]
 The Zak transform extends uniquely to
 a unitary operator ${\cal Z}: {\rm L}^2(\RR)
 \to {\cal H}$.
 \end{prop}
 
 The next lemma describes how the representation $\pi$ operates on the
 Zak transform side. It is easily verified on ${\cal Z}(C_c(\RR))$, and
 extends to ${\cal H}$ by density.
 \begin{prop}
  Let $\widehat{\pi}$ be the representation acting on ${\cal H}$,
  obtained by conjugating $\pi$ with ${\cal Z}$, i.e., $\widehat{\pi}
 (n,k,\overline{\ell}) = {\cal Z} \circ \pi(n,k,\overline{\ell})
 \circ {\cal Z}^*$. Then
 \begin{equation} \label{eqn:pi_hat}
 \widehat{\pi}(n,k,\overline{\ell})F (x,\omega) =
 e^{2 \pi i (\ell - nk)/L} e^{2 \pi i nx / L} F(x-k, \omega -
 n/L)~~.
 \end{equation}
 \end{prop}

 Now we can cite the Zak transform criterion for normalized
 tight Weyl-Heisenberg frames with integer oversampling. For a sketch of the
 proof confer \cite{Gr}, more details are contained in \cite{Da}.
 \begin{thm} \label{thm:ntf_Zak_crit}
  Let $f \in {\rm L}^2(\RR)$. Then $(M_{n/L} T_k f)_{n,k \in \ZZ}$ is
 a normalized tight frame of ${\rm L}^2(\RR)$ iff 
 \begin{equation} \label{eqn:ntf_Zak_crit}
 \sum_{i=0}^{L-1} | {\cal Z} f (x,\omega+i/L) |^2 = 1 ~~ \mbox{almost everywhere.}
 \end{equation}
 \end{thm}
 There exist more general versions of this criterion, which allow more
 complicated sets of time-frequency translations for the construction of the
 Gabor frames. While we have restricted our attention to the simple time-frequency
 lattice $\ZZ \times (1/L) \ZZ$ mostly for reasons of notational simplicity,
 the more general statements can be obtained employing suitable symplectic 
 automorphisms of the time-frequency plane. 

\section{Computing the Plancherel formula}

In this section we compute the Plancherel measure of $G $. 
Our calculations follow the recipe provided by Kleppner and Lipsman \cite{KlLi},
which is a natural extension of Mackey's procedure. Recall that the
Mackey machine allows to compute the dual of the group extension $G 
\supset N $ from the orbit space of the natural action of the quotient group
$G  / N $ and the duals of the associated fixed groups. For a detailed
account of the Mackey machine confer \cite{Fo}. In the following, we will
not explicitly distinguish between representations and their equivalence
classes. 

\noindent
{\bf Computing $\pmb{\widehat{G }}$:}\\
 We first note that since $G/N$ is finite, $N$ is regularly embedded in $G$,
 which is the chief technical requirement for the Mackey machine to run smoothly.
 Since $N$ is the direct product of three cyclic groups,
the character group $\widehat{N }$ is conveniently
parametrized by $[0,1[ \times [0,1[ \times \{ 0, 1, \ldots, L-1 \}$, by letting
\[
 \chi_{\omega_1,\omega_2,j}(nL,k,\overline{\ell}) = e^{2 \pi i (\omega_1 n + \omega_2 k + 
 j \ell /L)}~~.
\]
$G $ acts on $N $ by conjugation, which lifts to a natural action on 
$\widehat{N }$. Since
\[ (n,k,\overline{\ell}) (n'L,k',\overline{\ell'}) (n,k,\overline{\ell})^{-1} = (n'L,k',
 \overline{\ell'+k'n})~~, \]
we compute the dual action as
\[ (n,k,\overline{\ell}) . (\omega_1,\omega_2,j) = (\omega_1,\omega_2+jn/L- \lfloor \omega_2+jn/L \rfloor,j) ~~.\]
Here $\lfloor x \rfloor$ denotes the largest integer $\le x$.
Hence, defining 
\[ \Omega_j = [0,1[ \times [0,{\rm gcd}(j,L)/L[ \times \{ j \}~~,\]
 a measurable transversal of the orbits under the dual action is given by
$\Omega = \bigcup_{j=0}^{L-1} \Omega_j$.
Here $ {\rm gcd}(j,L)$ is the greatest common divisor of $j$ and $L$. The fact that
the subgroup of $\ZZ / L \ZZ$ generated by $\overline{j}$ coincides with the
subgroup generated by $\overline{{\rm gcd}(j,L)}$ accounts for this choice of transversal.
With the respect to the dual action, $(\omega_1,\omega_2,j) \in \Omega_j$ has
$N_j  = \{ ( n L/{\rm gcd}(j,L), k, \overline{\ell}) : k,n, \ell \in \ZZ \}$ as fixed group.
The associated {\bf little fixed group} is $N_j /N  \cong \ZZ / {\rm gcd}(j,L) \ZZ$.
For a convenient parametrisation of $\widehat{G }$ in terms of $\Omega$ and the
duals of the $N_j$ we need to establish the following lemma, which is verified
by straightforward calculation. This step is necessary because the extension $G 
\supset N $ is not a semidirect product. 
\begin{lemma}
 Let $(\omega_1,\omega_2,j) \in \Omega_j$ and $m \in \{0,1, \ldots, {\rm gcd}(j,L)-1 \}$.
 Then
 \[ \rho_{m,\omega_1,\omega_2,j}(n L/{\rm gcd}(j,L),k,\overline{\ell}) = e^{2 \pi i ((\omega_1 + m) n /{\rm gcd}(j,L)
 + \omega_2 k + j \ell/L)} ~~\]
 defines a character of $N_j$ with $\rho_{m,\omega_1,\omega_2,j}|_N = \chi_{\omega_1,\omega_2,j}$.
 Moreover, every irreducible representation of $N_j$ whose restriction to $N $ is a multiple of
 $\chi_{\omega_1,\omega_2,j}$ is equivalent to some $\rho_{m,\omega_1,\omega_2,j}$.
\end{lemma}
\begin{prf}
 The character property is verified by straightforward computation.
 The last statement is \cite[Proposition 6.40]{Fo}.
\end{prf}

Note that the additional parameter $m$ indexes the characters of the little fixed group $N_j / N $. 
Now the dual is obtained by inducing the $\rho_{m,\pmb{\omega},j}$, as follows from \cite[Theorems 6.38,6.39]{Fo}.
\begin{thm}
\label{thm:dual}
 Define, for $(\omega_1,\omega_2,j) \in \Omega$ and $m \in \{ 0, \ldots, {\rm gcd}(j,L)-1\}$
 the representation 
 \[ \sigma_{m,\omega_1,\omega_2,j} =  {\rm Ind}_{N_j}^{G } \rho_{m,\omega_1,\omega_2,j} ~~.\]
 If we let
 \[ \Sigma_j = \{ \sigma_{m,\omega_1,\omega_2,j} : (\omega_1,\omega_2,j) \in \Omega_j,
 m \in \{0,1, \ldots, {\rm gcd}(j,L)-1 \} ~~,\]
 then the dual of $G $ is the disjoint union
 \[
 \widehat{G } = \bigcup_{j =0}^{L-1} \Sigma_j
 \]
\end{thm}

\noindent
{\bf Computing Plancherel measure:}\\ As a preliminary remark we stipulate that the Haar measures
on all discrete groups $H$ occurring here are counting measures, normalized by $\mu_H(\{ e \})=1$.
This choice fixes the Plancherel measures uniquely, and implies in particular for all abelian
groups $H$ arising in the following that $\nu_H(\widehat{H})=1$.

 The reference for the following calculations is 
\cite{KlLi}, in particular \cite[II, Theorem 2.4]{KlLi}. Our computation follows the
proof to that result. Note that while the cited theorem only refers to the measure class,
the proof in fact gives the precise normalization. The starting point of the construction is
Theorem \ref{thm:dual}, which allows the identification
\[ \widehat{G } = \bigcup_{j=0}^{L-1} \Sigma_j= \bigcup_{j=0}^{L-1} \left( N_j / N  \right)^\wedge \times \Omega_j ~~.\]
On each of the $\Sigma_j$, Plancherel measure is a product measure: The
$ \left( N_j / N  \right)^\wedge$ carry the Plancherel measure of the finite quotient group,
which is simply counting measure weighted with $1/|N_j/N| = 1/{\rm gcd}(j,L)$.
For the missing parts, we decompose
Plancherel measure of $N$ on $\widehat{N }$ along orbits of the dual action.
This results in a measure on $\Omega \simeq \widehat{N}/G$, and the restrictions
to the $\Omega_j$ provide the second factors.
In order to explicitly compute these we note that
the Plancherel measure on $\widehat{N } \cong [0,1[ \times [0,1[ \times \{ 0, 1,
 \ldots, L-1 \}$ is $1/L$ times the product measure of Lebesgue measure on the first two factors and counting
measure on the third.
Since each orbit carries counting measure, the measure on the quotient is simply
Lebesgue measure on the transversal $ [0,1[ \times [0,{\rm gcd}(j,L)/L[$, for each $j$.
Thus we arrive at:

\begin{thm}
 \label{thm:Pl_meas}
 The Plancherel measure of $G $ is given by
 \begin{equation} \label{eqn:Pl_meas}
 d\nu_G(\sigma_{m,\omega_1,\omega_2,j}) = \frac{1}{L {\rm gcd}(j,L)} dm d\omega_1 d\omega_2 dj ~~.\end{equation}
 Here $d\omega_1$ and $d\omega_2$ are
 Lebesgue measure on the intervals $[0,1[$ and $[0,{\rm gcd}(j,L)/L[$, and $dm, dj$ are counting measure
 on $\{0,\ldots, {\rm gcd}(j,L) -1 \}$ and $\{ 0,\ldots, L-1\}$, respectively.
\end{thm}

As we will see in the next section, only the set $\Sigma_1$ will be of interest
for the Weyl-Heisenberg frame setting. Here the indexing somewhat simplifies: $N_1=N$, and 
$m$ can only take the value $0$. So we can identify $\Sigma_1$ with $\{ 0 \} \times [0,1[ \times [0,1/L[ \times \{1 \} $.

\section{Zak transform and Plancherel transform}

The aim in this section is to exhibit the representation $\widehat{\pi}$ obtained
by conjugating $\pi$ with the Zak transform as a direct integral of irreducibles.
This is done by taking a second look at (\ref{eqn:pi_hat}), which is a twisted
action by translations along $\ZZ \times (1/L) \ZZ$. Hence 
a decomposition of Lebesgue measure along cosets of $\ZZ \times (1/L) \ZZ$ gives
rise to a decomposition into representations acting on the cosets, and the
twisted action of the latter representations reveals them as induced representations.

To make this more precise, we let for $\pmb{\omega} \in [0,1[ \times [0,1/L[$ 
denote ${\cal O}_{\pmb{\omega}} = \pmb{\omega} + \ZZ \times (1/L) \ZZ$. The 
following lemma reveals the direct integral structure of $\widehat{\pi}$,
just by making suitable identifications.

\begin{lemma}
 Define for $\pmb{\omega} \in [0,1[ \times [0,1/L[$ the Hilbert space
\[ {\cal H}_{\pmb{\omega}} = \{ F : {\cal O}_{\pmb{\omega}} \to \CC : F \mbox{ fulfills }
 (\ref{eqn:per_cond}) \} ~~,\]
 with the norm defined by 
\begin{equation} \label{eqn:fibre_norm}
 \| F \|_{{\cal H}_{\pmb{\omega}}}^2 = \sum_{i=0}^{L-1} | F(\pmb{\omega} + (0,i/L)) |^2 .\end{equation}
 Let $\widehat{\pi}_{\pmb{\omega}}$ be the representation acting on ${\cal H}_{\pmb{\omega}}$
 by
 \[  \widehat{\pi}_{\pmb{\omega}} (k,n,\overline{\ell})F (\pmb{\gamma}) =
 e^{2 \pi i (\ell + nk)/L} e^{2 \pi i nx / L} F(\pmb{\gamma}-(k, n/L))~~.
 \]
 Then 
 \begin{equation} \label{eqn:dir_int_1}
 \widehat{\pi} \simeq \int_{[0,1[ \times [0,1/L[}^\oplus \widehat{\pi}_{\pmb{\omega}} ~d\pmb{\omega} ~~,
 \end{equation}
 via the map 
 \begin{equation} \label{eqn:intertwiner}
 F \mapsto \left( F |_{{\cal O}_{\pmb{\omega}}} \right)_{\pmb{\omega} \in [0,1[ \times [0,1/L[}
 \end{equation}
\end{lemma}

Strictly speaking, the intertwining operator (\ref{eqn:intertwiner}) is not well-defined
for arbitrary $F \in {\cal H}$, since the ${\cal O}_{\pmb{\omega}}$ are nullsets. However, 
the definition is rigorous for continuous $F$ and extends by density.
As a first glimpse of the connection between conditions (\ref{eqn:ntf_Zak_crit}) and 
(\ref{eqn:adm_cond_Pl}) note that the right-hand side of (\ref{eqn:ntf_Zak_crit})
can now be reformulated as
\[ \left\| \left( {\cal Z}f\right)|_{{\cal O}_{\pmb{\omega}}} \right\|_{{\cal H}_{\pmb{\omega}}} = 1, ~~
 \mbox{for almost every } \pmb{\omega} \in [0,1[ \times [0,1/L[~~.\] 
Hence the final step is to note that
(\ref{eqn:dir_int_1}) is in fact a decomposition into irreducibles:

\begin{lemma}
 If $\pmb{\omega} \in [0,1[ \times [0,1/L[$, then
 $\widehat{\pi}_{\pmb{\omega}} \simeq \sigma_{0,\pmb{\omega},1} \in \Sigma_1$. 
\end{lemma}

\begin{prf}
 We will use the imprimitivity theorem to show that $\widehat{\pi}$ is
 induced from a character of $N $. For this
 purpose consider the set $S = \{ \pmb{\omega} + (0,i/L) : i = 0,\ldots, L-1 \}$,
 with an action of $G $ on $S$ given by
 \[
   (n,k,\overline{\ell}) . \pmb{\gamma} = (\gamma_1,\gamma_2-n/L- \lfloor \gamma_2-n/L \rfloor)~~.
 \]
 The action is transitive with $N$ as associated stabilizer.
 To any subset $A \subset S$ we associate a projection operator $P_A$ on ${\cal H}_{\pmb{\omega}}$
 defined by pointwise multiplication with the characteristic function of 
 $A + \ZZ \times \ZZ$. It is then straightforward to check that $A \mapsto P_A$ is a
 projection-valued measure on $S$ satisfying
 \[ \widehat{\pi}_{\pmb{\omega}} (n,k,\overline{\ell}) P_A 
 \widehat{\pi}_{\pmb{\omega}} (n,k,\overline{\ell})^* = P_{(n,k,\ell).A} ~~.\]
 In other words, $A \mapsto P_A$ defines a transitive system of imprimitivity. Hence the
 imprimitivity theorem \cite[Theorem 6.31]{Fo} applies to show that $\widehat{\pi}_{\pmb{\omega}}
 \simeq {\rm Ind}_{N }^{G } \rho$ for a suitable representation $\rho$ of $N $. 
 Since the system of imprimitivity is based on a discrete set, we can follow the
 procedure outlined in \cite{Fo} immediately after Theorem 6.31, which identifies
 $\rho$ as the representation of $N$ acting on $P_{\{ \pmb{\omega} \}}({\cal H}_{\omega})$.
 For this purpose consider the function $F \in {\cal H}_{\pmb{\omega}}$ defined by 
 \[ F(\pmb{\omega}+(0,m/L)) = \delta_{m,0} ~~\mbox{for}~m=0,\ldots, L-1 ~~.\]
 Now the fact that
 \[ \widehat{\pi}_{\pmb{\omega}} (nL,k,\overline{\ell}) F = e^{2 \pi i \ell/L} e^{2 \pi i \omega_1 n}
 e^{2 \pi i \omega_2 k} F = \chi_{\omega_1,\omega_2,1} (nL,k,\overline{\ell}) F \]
 shows that 
 \[ \widehat{\pi}_{\pmb{\omega}} \simeq {\rm Ind}_N^G \chi_{\omega_1,\omega_2,1}
 = \sigma_{0,\omega_1,\omega_2,1,} ~~.\]
\end{prf}

By the last lemma and Mackey's theory, no two representations appearing 
in (\ref{eqn:dir_int_1}) are equivalent. Since $\pi$ is type-I, it follows that the
commuting algebra of $\pi$ is diagonal with respect to (\ref{eqn:dir_int_1}), hence abelian. But
this means that $\pi$ is multiplicity-free. 
Now a comparison of (\ref{eqn:dir_int_1}) with the Plancherel decomposition 
gives the desired result:
\begin{cor}
\label{cor:main}
 The Zak transform intertwines $\pi$ with a direct integral
 \[ \widehat{\pi} \simeq \int_{\Sigma_1}^\oplus \sigma ~d\nu_{G }(\sigma) ~~.\]
 The criterion (\ref{eqn:ntf_Zak_crit}) is an immediate consequence
 of Lemma \ref{lem:ntf_eq_admissible} and the admissibility condition (\ref{eqn:adm_cond_Pl}). 
\end{cor}

\section{Concluding remarks}

 The connection between Weyl-Heisenberg criteria and direct integrals
has already been addressed by other authors, more or less explicitly,
see for instance \cite{RoSh,Wo}. However, we are not aware of any
previous reference to nonabelian Plancherel theory in this context.
The case of rational oversampling, which can also be
dealt with using the Zak transform, does not seem to fit into the 
Plancherel transform setting as neatly as the integer oversampling case does
according to Corollary \ref{cor:main}. 

Even more complicated is the case of irrational oversampling, which
amounts to replacing $L$ in the above definitions by 
an irrational $\alpha>1$. Again it is simple to establish that
normalized tight frame conditions are equivalent to admissibility with
respect to a certain representation of a suitably defined group $G$.
The intriguing fact about this representation is that while it
can be shown to be not of type I, there exists a characterization of the 
admissible vectors, due to Ron and Shen \cite[Corollary 2.19]{RoSh},
which is very similar to the admissibility condition formulated in
\cite{Fu}, after making suitable identifications.
In a sense, the Plancherel measure is replaced by a family
of measures, each effecting a decomposition into irreducibles.
Now admissible vectors have to fulfill conditions with respect
to each of these measures which are entirely analogous to the
admissibility condition involving Plancherel measure in the
type-I case. And conversely,
the joint admissibility conditions are also sufficient.
A better understanding of this example should provide some orientation
for dealing with more general groups. A more detailed exposition
of the connections between irrational oversampling and non type-I
admissibility criteria will be given elsewhere.


\begin{thebibliography}{99}
 \bibitem{BeTa}{D. Bernier and K. Taylor, {\em Wavelets from 
  square-integrable representations.} SIAM J. Math. Anal. {\bf 27} (1996), 
  594-608.}
\bibitem{Da} {I. Daubechies, {\em The wavelet transform, time-frequency
 localization and signal analysis.} 
 IEEE Trans. Inform. Theory {\bf 34} (1988), 961-1005.}
\bibitem{Di} {J. Dixmier, {\em $C^{\ast}$-Algebras.}
 North Holland, Amsterdam, 1977.}
\bibitem{Fo} {G.B. Folland, {\em A Course in Abstract Harmonic Analysis.}
 CRC Press, Boca Raton, 1995.}
\bibitem{Fu96} {H. F\"uhr, {\em Wavelet frames and admissibility in higher
 dimensions.} J. Math. Phys. {\bf 37} (1996), 6353-6366.}
\bibitem{FuMa} {H. F\"uhr and M. Mayer, {\em Continuous wavelet transforms
 from semidirect products: Cyclic representations and Plancherel measure.}
 To appear in J. Fourier Anal. Appl. Electronically available as
 \underline{\sf math-ph/0102002}.}
\bibitem{Fu} {H. F\"uhr, {\em Admissible vectors for the regular
 representation.} To appear in Proc AMS. Electronically available as
 \underline{\sf math-ph/0010051}.}
 \bibitem{Gr}{K. Gr\"ochenig, {\em Foundations of Time-Frequency Analysis.}
  Birkh\"auser, Boston, 2001. }
 \bibitem{Ka}{E. Kaniuth, {\em Der Typ der regul\"aren Darstellung
 diskreter Gruppen.} Math. Ann. {\bf 182} (1969), 334--339.} 
 \bibitem{KlLi} {A. Kleppner and R.L. Lipsman, {\em The Plancherel formula
  for group extensions, I and II.} Ann.Sci.Ecole Norm.Sup. {\bf 5} (1972),
  459-516; ibid. {\bf 6} (1973), 103-132.}
  \bibitem{LWWW}{R.S. Laugesen, N. Weaver, G. Weiss and E.N. Wilson,
 {\em Continuous wavelets associated with a general class of 
 admissible groups and their characterization.} J. Geom. Anal., to appear.}
\bibitem{RoSh}{A. Ron and Z. Shen, {\em Weyl-Heisenberg frames and Riesz bases in
 ${\rm L}^2(\RR^d)$}, Duke Math. J. {\bf 89 (2)} (1997), 237-282.}
\bibitem{Th}{E. Thoma, {\em Eine Charakterisierung diskreter Gruppen vom Typ I}, Invent. Math. {\bf 6} (1968), 
 190-196.}
\bibitem{Wo}{P. Wojdy{\l}{\l}o,{\em Gabor and Wavelet Frames}, Thesis, Warsaw University, 2000.}

\end{thebibliography}
\end{document}